
\documentstyle{amsppt}
\magnification=\magstep 1


\define\codim{\operatorname{codim}}
\define\rank{\operatorname{rank}}
\define\Sym{\operatorname{Sym}}
\define\Adams{\operatorname{Adams}}
\define\Ann{\operatorname{Ann}}
\NoBlackBoxes

\topmatter
\title
Residual Intersections and Some Applications
\endtitle{}
\author Xian Wu
\endauthor{}
\address
Department of Mathematics,
University of South Carolina,
Columbia,
SC 29208
\endaddress{}
\email
wux{\@}milo.math.scarolina.edu
\endemail{}

\abstract
We give a new residual intersection decomposition
for the refined intersection products of Fulton-MacPherson.
Our formula refines the celebrated residual intersection formula
of Fulton, Kleiman, Laksov, and MacPherson.
The new decomposition is more likely to be compatible with
the canonical decomposition of the intersection products
and each term in the decomposition thus has simple geometric meaning.
Our study is motivated by its applications to some geometric problems.
In particular,
we use the decomposition
to find the distribution of limiting linear subspaces
in degenerations of hypersurfaces.
A family of identities for characteristic classes
of vector bundles is also obtained as another consequence.
\endabstract{}


{}

\endtopmatter{}
\document

\heading 1. Introduction
\endheading{}

Given a closed regular embedding of codimension $d$
$$
i\ : \ X \to Y
$$
and a morphism
$$
f\ : \ V \to Y
$$
with $V$ a purely $k$-dimensional scheme,
the fundamental construction of Fulton-MacPherson
[F] defines the refined intersection product
$$
X\cdot V \ \in A_{k-d}(W),\quad W=f^{-1}(X),
$$
where $A_{k-d}(W)$ is the $(k-d)$-th Chow group of $W$.
Let $N$ be the pull-back of the normal bundle $N(X,Y)$ of $X$ in $Y$
to $W$ and $c(N)$ be its Chern class.
Then $X\cdot V$ can be expressed in terms
of $c(N)$ and the Segre class $s(W,V)$ of $W$ in $V$ by
$$
X\cdot V=\{c(N)\cap s(W,V)\}_{k-d}.
$$
Furthermore,
there is a canonical decomposition
$$
X\cdot V = \sum_{j}\alpha_j,
\tag 1.1
$$
where $\alpha_j$ are classes supported on the so-called
distinguished varieties $Z_j$ of $X\cdot V$.
It is well-known that every irreducible component of $W$
is a distinguished variety of $X\cdot V$.
The canonical decomposition gives all the important
information about $X\cdot V$.
However,
it is in general very difficult to find such a decomposition.

There are many other ways to decompose
$X\cdot V$
and in many applications a less comprehensive decomposition
will be enough.
Given a closed subscheme $Z$ of $W$,
a basic problem of residual intersections
is to decompose $X\cdot V$ into the sum of a class
supported on $Z$ and a class supported on a ``residual''
set of $W$ with respect to $Z$.
The celebrated residual intersection formula of
Fulton, Kleiman, Laksov, and MacPherson
gives such a decomposition
$$
X\cdot V=\{c(N)\cap s(Z,V)\}_{k-d}+\Bbb R =: \Bbb M_Z+\Bbb R,
\tag 1.2
$$
where $s(Z,V)$ is the Segre class of $Z$ in $V$.
The construction of the residual intersection class $\Bbb R$
is more involved and one can refer to [F] for details.
This decomposition is much easier to compute than the
canonical decomposition and
the formula has many important applications.
Please see Chapter 9 of [F] for further references
on the extensive literature and
numerous results on residual intersections.

If $Z$ is a connected component of $W$,
then $\Bbb M_Z$ in (1.2) is nothing but the term (or the sum of terms)
in the canonical
decomposition (1.1) supported on $Z$ and $\Bbb R$
is the sum of all other terms in the right-hand side of (1.1).
However,
those two decompositions are not compatible in general.
For example,
assume that $W$ is the union of two irreducible components $Z$ and $R$
and there are no contributions to $X\cdot V$ other than those from
these two components.
In this case,
both decompositions have two terms which are
supported on the same varieties $Z$ and $R$, respectively.
However,
they are generally not the same unless $Z\cap R$ has a dimension
less than $k-d$.
In fact,
there seems to be no good coordination
between two decompositions in general.
Consider the case when $N$ is generated by its sections.
Then it is well-known that
all the terms in (1.1) must be represented by non-negative cycles.
Furthermore,
they can be interpreted geometrically as limits
by using dynamic intersections and deformation theory.
On the other hand,
as we will see in Section 4 (Example 4.14),
both terms in the right-hand side of
(1.2) can be negative (of course, not at the same time).
Therefore,
(1.2) cannot be obtained by shuffling terms of (1.1)
and it is hard to interpret the formula geometrically.
One cause of this is easy to see,
since two terms in the right-hand side of
(1.2) are not equally weighted.
The formula for $\Bbb M_Z$
uses no information about the residual set $R$
while construction of $\Bbb R$ depends
on structures of both $Z$ and $R$.
In a sense,
(1.2) packs all the complication
into $\Bbb R$ to make $\Bbb M_Z$ as simple as possible.
While this is very useful in many applications,
it also makes the formula less favorable in other situations.

In this paper,
we propose another way to decompose $X\cdot V$
into the sum of a class supported on $Z$ and a class supported
on a ``residual'' subscheme of $W$ with respect to $Z$.
Our decomposition is somehow ``symmetric''.
This makes it more likely to be compatible with the canonical
decomposition and the decomposition can thus
be used to find the equivalence of $Z$
for $X\cdot V$.
In fact,
in the case that (1.1) has only two terms as
mentioned in the last paragraph,
our decomposition gives the actual canonical decomposition.
In particular,
it has nice geometric meaning as the distribution of the limits.
Of course,
there is a price to be paid for this.
Mainly,
both terms in the decomposition are now equally complicated.
However,
it is still manageable and can be
computed (with the help of computers) in many cases.
We are particularly encouraged by its applications for some
geometric problems that initially motivated this study.

This paper is organized as follows.
In Section 2,
one can find basic formulas for our new decomposition
(Theorem 2.4 and Corollary 2.10).
The idea behind is very simple.
If $W$ is the sum of two divisors in $V$,
then there is a very natural way to decompose $X\cdot V$ symmetrically
into the sum of two classes that are supported on the divisors.
In general,
if a closed subscheme $Z$ of $W$ is given,
one can define a residual scheme $R(Z)$ to $Z$ in $W$
such that
$$
W= Z \cup R(Z).
\tag 1.3
$$
We then decompose $X\cdot V$ into the sum of a class supported
on $Z$ and a class supported on $R(Z)$
by reducing to the case of divisors using blow-ups.
Our formulas are given in such a way so that they can be
easily compared with the standard residual intersection
formula (1.2) of Fulton, Kleiman, Laksov, and MacPherson.
In particular,
we see that the new decomposition refines the standard one
in a sense that will be made clear in Section 2.

A closer look of the new decomposition is given
in Section 3.
In general,
$R(R(Z))$ may not be equal to $Z$.
Therefore,
(1.3) and hence the decomposition of $X\cdot V$ given in Section 2
will not be symmetric in general.
Even in the case that (1.3) is symmetric,
our decomposition of $X \cdot V$ still depends on the order of
two blow-ups needed in the process (see Example 3.1).
In fact,
as in the standard decomposition (1.2),
it can still happen that neither term in
the new one is related to the scheme
structure of $R(Z)$.
In spite of all of that,
the new intersection decomposition
does behave well
for many interesting decompositions
$$
W= Z_1 \cup Z_2
$$
of $W$
(Theorem 3.4).
In particular,
we obtain some new formulas which express our decomposition
of $X\cdot V$ in a truly symmetric form in terms that are related to
the scheme structures of $Z_1$ and $Z_2$
under some modest hypotheses (Theorem 3.6 and Corollary 3.16).

The main motivation that inspired
our search for a refined residual intersection
decomposition is its geometric meaning.
This and some applications of the new decomposition
are studied in Section 4.
If we interpret the refined intersection product $X\cdot V$
as the class of the limits derived from deformation theory,
then our decomposition tells us how the limits
are distributed in different components.
In particular,
we use it to study degenerations of hypersurfaces in $\Bbb P^n$
and their limiting $\Bbb P^r$'s (Proposition 4.1).
This not only provides a much simpler and intrinsic way to recapture
some results of [W1], [W2], and [W3],
but also yields new results which would be very difficult
to obtain from deformation theory directly.
Actual calculations of some examples can also be found there.
As another application,
we obtain a family of identities for characteristic classes
of vector bundles.

\example{Acknowledgement}
The author is indebted to W. Fulton for valuable suggestions
and advice.
The author also would like to
thank P. Aluffi, A. Kustin, and M. Miller for helpful discussions.
Finally,
the author wants to thank the referee for a suggestion about
a generalization of Corollary 4.15.
\endexample{}

\example{Blanket Conventions}
To simplify the notation,
we will make the following conventions for this paper.
Hopefully,
the meaning will be clear from the content.

(1) A class in the Chow group of
a subspace will often be identified
with its image in the Chow group of a bigger space
without warning.

(2) The same notation will be used to denote a vector bundle
and its restriction to any subspace.

(3) Similarly,
the same notation will be used to denote a morphism and its
restriction to any subspace.

(4) Finally,
the same notation will be used to denote a divisor,
the corresponding line bundle,
and the first Chern class of the bundle.
\endexample{}

\smallpagebreak
\heading 2. A refined residual intersection formula
\endheading{}
\smallpagebreak

Let us first recall the basic construction of refined intersection
products of Fulton-MacPherson from [F].
Consider the fiber square
$$
\CD
W @>>> V\\
@VfVV @VVfV\\
X @>>i> Y,
\endCD{}
\tag 2.1
$$
where
$$
i\ : \ X \to Y
$$
is a closed regular embedding of codimension $d$,
$$
f\ : \ V \to Y
$$
is a morphism
with $V$ a purely $k$-dimensional scheme,
and $W$ is the inverse image scheme $f^{-1}(X)$.
{}From this,
the fundamental construction of Fulton-MacPherson
defines the refined intersection product
$X\cdot V$ in the $(k-d)$-th Chow group $A_{k-d}(W)$
of $W$ by embedding the normal cone $C_WV$ of $W$ in $V$
into the vector bundle $N$,
where $N$ is the pull-back of the normal bundle $N(X,Y)$ of $X$ in $Y$
to $W$.
Let $c(N)$ be the Chern class of $N$
and $s(W,V)$ be the Segre class of $W$ in $V$.
Then $X\cdot V$ can be expressed as
$$
X\cdot V=\{c(N)\cap s(W,V)\}_{k-d}.
\tag 2.2
$$
Given a closed subscheme $Z$ of $W$,
our goal in this section is,
for a suitable decomposition
$$
W=Z \cup Z'
$$
of $W$,
to decompose
$X\cdot V$ in a natural way
as the sum of two classes $\Bbb R_{Z}$ and $\Bbb R_{Z'}$ that
are supported on $Z$ and $Z'$,
respectively.
To simplify the notation,
we will assume that $W$ has at least codimension 1 in $V$.
The case of $W=V$ is simple since we then have
$$
s(W,V)= [V] \quad \text {and} \quad X \cdot V = c_{top}(N)\cap [V].
$$

Consider first the case that $Z$ is a divisor in $V$.
Let $R(Z)$ be the residual scheme to $Z$ in $W$.
To be more precise,
in general,
given a closed subscheme $Z$ of $W$,
we define the residual subscheme $R(Z)$ to $Z$ in $W$
with respect to $V$
to be the subscheme of $V$ defined by the ideal sheaf
$$
\Cal I_{R(Z)}=\Ann(\Cal I_Z/\Cal I_W),
$$
where $\Cal I_Z$ and $\Cal I_W$ are ideal sheaves of $Z$ and $W$ in $V$,
respectively.
We hence have set-theoretically
$$
W=Z \cup R(Z).
\tag 2.3
$$
Moreover,
in the case that $Z$ is a divisor in $V$,
(2.3) gives actually a scheme-theoretic decomposition of $W$
in the sense that
$$
\Cal I_W = \Cal I_Z \cdot \Cal I_{R(Z)}.
$$
For a given vector bundle $E$,
the standard notation $c_i(E)$ for the $i$-th Chern class of $E$
and $s_i(E)$ for the $i$-th Segre class of $E$
will be used.
Recall also
that we will use the same notation to denote a divisor and
its corresponding line bundle.
We are now ready to propose the following residual
intersection decomposition of $X\cdot V$.

\proclaim{Theorem 2.4}
Consider the following expansion of the fiber square (2.1)
$$
\CD
@. R @.\\
@. @VVV @.\\
D @>>> W @>>> V\\
@. @VfVV @VVfV\\
@. X @>>i> Y,
\endCD{}
\tag 2.5
$$
where the subscheme $D$ of $W$ is a divisor in $V$
and $R$ is the residual scheme to $D$ in $W$.
Let
$\Bbb R_{D}$ and $\Bbb R_{R}$ be two
classes supported on
$D$ and $R$ defined by
$$
\Bbb R_D=\{c(N)\cap s(D,V)\}_{k-d}
+\sum_{i=0}^{d-2}\sum_{j=1}^{d-1-i}{d-1-i \choose j}
c_i(N)s_{d-i-j}(D)\cap s_{k-j}(R,V)
\tag 2.6
$$
and
$$
\Bbb R_R=\{c(N)\cap s(R,V)\}_{k-d}
+\sum_{i=0}^{d-2}\sum_{j=1}^{d-1-i}{d-1-i \choose j}
c_i(N)s_j(D)\cap s_{k-d+i+j}(R,V),
\tag 2.7
$$
respectively.
Then
the refined intersection product $X\cdot V$ can be decomposed into
$$
X\cdot V = \Bbb R_{D} + \Bbb R_{R}.
\tag 2.8
$$
\endproclaim{}

\example{Remark}
There are many different ways to express $\Bbb R_D$ and $\Bbb R_R$
and they may serve some particular needs better than the forms given above.
In the proof given below,
one can find some different formulas of $\Bbb R_D$ and $\Bbb R_R$
that are less messy.
We choose the formulas (2.6) and (2.7) not only
because they are very explicit
but also to compare our decomposition with the
standard residual intersection decomposition (1.2) as given in [F].
\endexample{}

\demo{Proof of Theorem 2.4}
By our assumption on $W$,
$R$ is a subscheme of codimension at least 1 in $V$.
Consider the blow-up of $V$ along $R$
$$
b\: \ \widetilde V \to V
$$
Let
$$
\widetilde W = b^{-1}(W),
\quad \widetilde D=b^{-1}(D)=b^*(D),
\quad \text {and} \quad \widetilde R=b^{-1}(R)
$$
be the exceptional divisor.
Consider the following fiber diagram:
$$
\CD
\widetilde W @>>> \widetilde V\\
@VbVV @VVbV\\
W @>>> V\\
@VfVV @VVfV\\
X @>>i> Y,
\endCD{}
$$
where the low half of the diagram is just the fiber square (2.1).
By the push-forward formula [Theorem 6.2, F],
we have
$$
b_*(X \cdot \widetilde V)
=b_*(i^!([\widetilde V]))
=i^!(b_*([\widetilde V]))
=i^!([V])
=X \cdot V,
$$
where $i^!$ is the refined Gysin homomorphism defined
from refined intersections.
Notice that $\widetilde W$ is now a divisor in $\widetilde V$
and there is a very natural way to decompose
$X\cdot \widetilde V$ symmetrically with respect to two components
$\widetilde D$
and
$\widetilde R$
of
$\widetilde W$.
In fact,
$$
\aligned
X \cdot V
&=b_*(X \cdot \widetilde V)\\
&=b_*(\{c(b^*N)\cap s(\widetilde D + \widetilde R,
\widetilde V)\}_{k-d})\\
&=b_*(\sum_{i=0}^{d-1}c_i(b^*N)
\cap s_{k-d+i}(\widetilde D + \widetilde R, \widetilde V))\\
(2.9) \quad &=b_*(\sum_{i=0}^{d-1}c_i(b^*N)
(-\widetilde D - \widetilde R)^{d-1-i}\cap
([\widetilde D + \widetilde R]))\\
&=\sum_{i=0}^{d-1}c_i(N)\cap
b_*((-\widetilde D -\widetilde R)^{d-1-i}\cap [\widetilde D])
+\sum_{i=0}^{d-1}c_i(N)\cap
b_*((-\widetilde D -\widetilde R)^{d-1-i}\cap [\widetilde R])\\
&=:\Bbb R_D + \Bbb R_R.
\endaligned{}
$$

Before going further,
let us point out that the formulas above are already
in computable forms in many cases.
They are truly symmetric with respect to
$\widetilde D$
and
$\widetilde R$.
The drawback is that they are expressed indirectly by push-forwards.

We will now exam each term more closely.
Let us start with $\Bbb R_R$.
$$
\aligned
\Bbb R_R&=\sum_{i=0}^{d-1}c_i(N)\cap
b_*((-\widetilde D -\widetilde R)^{d-1-i}\cap [\widetilde R])\\
&=\sum_{i=0}^{d-1}c_i(N)\cap
b_*(\sum_{j=0}^{d-1-i}{d-1-i \choose j}
(-\widetilde D)^j (-\widetilde R)^{d-1-i-j}\cap [\widetilde R])\\
&=\sum_{i=0}^{d-1}c_i(N)\cap
b_*(\sum_{j=0}^{d-1-i}{d-1-i \choose j}
(-b^*D)^j\cap (-\widetilde R)^{d-1-i-j}\cap [b^* R])\\
&=\sum_{i=0}^{d-1}\sum_{j=0}^{d-1-i}{d-1-i \choose j}
c_i(N)s_j(D)\cap s_{k-d+i+j}(R,V)\\
&=\sum_{i=0}^{d-1}c_i(N)\cap s_{k-d+i}(R,V)
+\sum_{i=0}^{d-2}\sum_{j=1}^{d-1-i}{d-1-i \choose j}
c_i(N)s_j(D)\cap s_{k-d+i+j}(R,V)\\
&=\{c(N)\cap s(R,V)\}_{k-d}
+\sum_{i=0}^{d-2}\sum_{j=1}^{d-1-i}{d-1-i \choose j}
c_i(N)s_j(D)\cap s_{k-d+i+j}(R,V).
\endaligned{}
$$
Similarly,
we can get the formula for $\Bbb R_D$ as follows.
$$
\aligned
\Bbb R_D&=\sum_{i=0}^{d-1}c_i(N)\cap
b_*((-\widetilde D -\widetilde R)^{d-1-i}\cap [\widetilde D])\\
&=\sum_{i=0}^{d-1}c_i(N)\cap
b_*(\sum_{j=0}^{d-1-i}{d-1-i \choose j}
(-\widetilde R)^j(-\widetilde D)^{d-1-i-j}\cap [\widetilde D])\\
&=\sum_{i=0}^{d-1}c_i(N)\cap s_{k-d+i}(D,V)
+\sum_{i=0}^{d-2}\sum_{j=1}^{d-1-i}{d-1-i \choose j}
c_i(N)b_*((-\widetilde R)^{j-1}(-\widetilde D)^{d-i-j}\cap [\widetilde R])\\
&=\{c(N)\cap s(D,V)\}_{k-d}
+\sum_{i=0}^{d-2}\sum_{j=1}^{d-1-i}{d-1-i \choose j}
c_i(N)s_{d-i-j}(D)\cap s_{k-j}(R,V).
\endaligned{}
$$
This completes our proof.
\qed
\enddemo{}

In general,
if a closed subscheme $Z$ of $W$ is given,
we can blow up $V$ along $Z$ to reduce to the case of
a divisor in the blow-up $\widetilde V$.
We then obtain
our decomposition by pushing forward the decomposition of
$X\cdot \widetilde V$ as stated in Theorem 2.4.

\proclaim{Corollary 2.10}
Let $X \cdot V$ be the refined intersection product
defined from the fiber square (2.1).
Given a closed subscheme $Z$ of $W$,
let
$$
\pi\ : \
\widetilde V \to V
$$
be the blow-up of $V$ along $Z$
and $D$ be the exceptional divisor.
Furthermore,
let $\widetilde R$ be the residual scheme to $D$ in $\widetilde V$
and $R$ be the image of $\widetilde R$ under $\pi$.
Then
$X\cdot V$ can be decomposed into
$$
X\cdot V = \Bbb R_{Z} + \Bbb R_{R},
\tag 2.11
$$
where
$\Bbb R_{Z}$ and $\Bbb R_{R}$ be two
classes supported on
$Z$ and $R$ defined by
$$
\Bbb R_Z= \Bbb M_Z + \Bbb A_Z,
\tag 2.12
$$
$$
\Bbb M_Z = \{c(N)\cap s(Z,V)\}_{k-d}
=\sum_{i=0}^dc_i(N) \cap s_{k-d+i}(Z,V),
\tag 2.13
$$
$$
\Bbb A_Z = \sum_{i=0}^{d-2}\sum_{j=1}^{d-1-i}{d-1-i \choose j}
c_i(N)\cap \pi_*(s_{d-i-j}(D)\cap s_{k-j}(\widetilde R,\widetilde V))
\tag 2.14
$$
and
$$
\Bbb R_R= \Bbb M_R + \Bbb A_R,
\tag 2.15
$$
$$
\Bbb M_R = \{c(N)\cap \pi_*s(\widetilde R,\widetilde V)\}_{k-d}
=\sum_{i=0}^{d}c_i(N)\cap \pi_*s_{k-d+i}(\widetilde R,\widetilde V),
\tag 2.16
$$
$$
\Bbb A_R=\sum_{i=0}^{d-2}\sum_{j=1}^{d-1-i}{d-1-i \choose j}
c_i(N)\cap \pi_*(s_j(D)\cap s_{k-d+i+j}(\widetilde R,\widetilde V)),
\tag 2.17
$$
respectively.
\endproclaim{}

\example{Remark}
Let $R(Z)$ be the residual scheme to $Z$ in $W$ as defined earlier.
Then $R(Z)$ may not be equal to $R$ even as a set.
However,
$R$ is always contained in $R(Z)$ and
$\Bbb R_R$ can hence be considered as a class supported on $R(Z)$
as well.
A more detailed study on this matter will be conducted
in the next section.
\endexample{}

In closing of this section,
we observe that each term in our decomposition
as stated in Corollary 2.10 is expressed in two terms itself.
In particular,
we see that $\Bbb R_{Z}$ contains two parts.
One part $\Bbb M_Z$ is nothing but the very term that appears in
the standard residual intersection formula (1.2)
and we will call it the main term.
The other more complicated one $\Bbb A_Z$
is a class supported on
$Z\cap R$ and we will call it the adjunct term.
While this term is packed in the residual intersection class $\Bbb R$
in the standard
residual decomposition (1.2),
it is now squeezed out in the new decomposition
as we try to somehow even up the distribution.
In this sense,
the new decomposition gives a refinement of the standard decomposition
(1.2).

\smallpagebreak
\heading 3. The decomposition in symmetric forms
\endheading{}
\smallpagebreak

The idea behind the residual decomposition given in the last section
is very simple.
We see that if $W$ is the sum of two divisors of $V$,
then there is a very natural way to decompose $X\cdot V$
into the sum of two classes that are supported on the divisors.
The decomposition in general cases is then obtained
by reducing to the case of divisors using blow-ups.
In general,
we need two blow-ups to do that and the order matters.
In other words,
although the new decomposition has a more balanced distribution
between its terms,
it is still not symmetric with respect to
the residual decomposition (2.3)
$$
W= Z \cup R(Z)
$$
of $W$.
This is not surprising,
of course,
since the above decomposition of $W$ itself is not symmetric.
We do not in general have that
$R(R(Z))$ is equal to $Z$.
However,
that is not the main obstacle.
In fact,
even if the decomposition of $W$ above is symmetric,
our decomposition of $X\cdot V$ may still depends on which component
to be chosen first.
Indeed,
let $\pi$ be the blow-up of $V$ along $Z$ and $D$ be the exceptional
divisor.
In general,
the image $\pi(R(D))$ is not equal to $R(Z)$ even as a set,
where $R(D)$ is the residual scheme to $D$ in $\pi^{-1}(W)$.
Therefore,
the second term $\Bbb R_{R}$ in our decomposition is not expressed
in terms of (and may indeed not be related to)
the scheme structure of $R(Z)$.
The following simple example illustrates that.

\example{Example 3.1}
Let $A$ and $B$ be lines in $\Bbb P^2$ defined by
$$
x=0 \quad \text {and} \quad
y=0,
$$
respectively.
Consider the following fiber square
$$
\CD
W @>>> V=\Bbb P^2\\
@VfVV @VVfV\\
X=2A \times 2B @>>i> Y=\Bbb P^2 \times \Bbb P^2,
\endCD{}
\tag 3.2
$$
where $f$ is the diagonal embedding.
Then
$$
\Cal I_W = \{x^2, y^2\} \quad \text {and} \quad
X \cdot V = 4p,
$$
where
$$
p=A \cap B \quad \text {and} \quad
\Cal I_p = \{x,y\}.
$$
It is easy to see that the residual scheme
$R(p)$ to $p$ in $W$ is defined by
$$
\Cal I_{R(p)} = \{x^2, y^2, xy\} = \Cal I_p^2.
$$
Furthermore,
we do have
$$
R(R(p)) = p.
$$
Therefore,
if we take
$$
Z_1= p \quad \text {and} \quad
Z_2 = R(p)
$$
then we have a ``symmetric'' decomposition
$$
W = Z_1 \cup Z_2
\tag 3.3
$$
of $W$.
(Notice that the multiplicity of $W$, $Z_1$, and $Z_2$
are 4, 1, and 3,
respectively.)
However,
our decomposition of $X\cdot V$
into the sum of two classes that are supported on
$Z_1$ and $Z_2$ will depend on
which $Z_i$ to be blown up first.
Let $\pi_i$ be the blow-ups of $W$ along $Z_i$
and $D_i$ be the exceptional divisors.
Then it is easy to see that
the residual scheme to $D_1$ in $\pi_1^{-1}(W)$ is $D_1$ itself and
the residual scheme to $D_2$ in $\pi_2^{-1}(W)$ is the empty set.
The corresponding decompositions are therefore
$$
X\cdot V = 2p + 2p
\quad \text {and} \quad
X\cdot V = 4p + 0,
$$
respectively.
(As a comparison,
the corresponding decompositions using the standard
residual intersection formula (1.2)
are
$$
X\cdot V = p + 3p
\quad \text {and} \quad
X\cdot V = 4p + 0,
$$
respectively.)
\endexample{}

While the examples such as one given above
make things more interesting,
they also suggest that the decomposition for $X\cdot V$
might behave ``better'' if
the decomposition of $W$ itself is ``nicer''.
As it turns out,
this is indeed the case
and the decompositions of $W$ in terms of residual schemes
are misleading.
In the following,
we will show that,
for some reasonable decompositions of $W$
as the union of two subschemes
$Z_1$ and $Z_2$,
our decomposition of $X\cdot V$ as given in Section 2
is actually independent of which component $Z_i$ to be chosen
first.
We will derive our formulas in truly
symmetric forms in terms
that are related to the scheme structures
of $Z_1$ and $Z_2$.

In many applications,
one is mainly interested in finding the contributions to $X\cdot V$
from certain components of $W$.
In this case,
the most direct way to reduce to the case of divisors
is to blow up $V$ along $W$.
In fact,
this is basically how $X\cdot V$ is defined.
This simple observation motivates the following theorem.

\proclaim{Theorem 3.4}
Let
$$
\pi\: \ \widetilde V \to V
$$
be the blow-up of $V$ along $W$.
Assume that a decomposition of
$$
W= Z_1 \cup Z_2
$$
is given such that the exceptional divisor
$$
E
=\pi^{-1}(W)
=\widetilde E_1 + \widetilde E_2,
$$
where
$\widetilde E_1$
and $\widetilde E_2$
are inverse image schemes of $Z_1$ and $Z_2$,
respectively,
then our refined
residual decomposition of $X \cdot V$ as defined in Section 2
is independent of the choice of which $Z_i$ to be considered first.
\endproclaim{}

In fact,
we will show something stronger.
Let
$$
\pi_1\: \ V_1 \to V
$$
be the blow-up of $V$ along $Z_1$ and
$$
E_1=\pi_1^{-1}(Z_1)
$$
be the exceptional divisor.
Furthermore,
let
$$
b_1\: \ V_1' \to V_1,
$$
be the blow-up of $V_1$ along $R_1$,
where $R_1$ is the residual scheme to $E_1$ in $\pi_1^{-1}(W)$,
and
$$
b_1^{-1}(\pi_1^{-1}(W))
=b_1^{-1}(E_1\cup R_1)
= E_1' + R_1',
$$
where $E_1'$ is the inverse image of $E_1$
under $b_1$ and $R_1'$ is the exceptional divisor.
Similarly,
we define blow-ups
$$
\pi_2\: \ V_2 \to V
\quad \text {and} \quad
b_2\: \ V_2' \to V_2
$$
and $E_2$, $R_2$, $E_2'$, and $R_2'$,
respectively.
Putting all those blow-ups together,
we consider the following fiber diagram:
$$
\CD
E_1'+R_1' @>>> V_1'
@.
@.\quad \quad \quad \quad\quad \quad \quad \quad
@.
E_2'+R_2' @>>> V_2'\\
@Vb_1VV @VVb_1V
@.@.@Vb_2VV @VVb_2V\\
E_1\cup R_1 @>>> V_1
@.
@.\quad \quad \quad \quad\quad \quad \quad \quad
@.
E_2\cup R_2 @>>> V_2
\endCD{}
$$
$$
\CD
\widetilde E_1 + \widetilde E_2 @>>> \widetilde V\\
@V\pi VV @VV\pi V \\
Z_1\cup Z_2 @>>> V\\
@VfVV @VVfV\\
X @>>i> Y
\endCD{}
\tag 3.5
$$
where
$p$ and $q$ are the unique
morphisms determined by the universal property of
blow-ups for $\pi$.
{}From the diagram above,
there are three different ways to decompose $X\cdot V$
into the sum of two classes supported on $Z_1$ and $Z_2$
corresponding to three push-forwards of the natural
decompositions of $X \cdot V_1'$,
$X \cdot \widetilde V$,
and $X \cdot V_2'$
as defined in Section 2.
We can now state the following stronger version of Theorem 3.4

\proclaim{Theorem 3.6}
All three procedures as described above give rise to
the same decomposition of $X\cdot V$.
In particular,
the decomposition of
$X\cdot V$ as defined in Section 2
can be expressed in the following symmetric form
$$
X\cdot V = \Bbb R_{Z_1} + \Bbb R_{Z_2}
\tag 3.7
$$
and $\Bbb R_{Z_1}$ and $\Bbb R_{Z_2}$ are
the classes supported on
$Z_1$ and $Z_2$ defined by
$$
\Bbb R_{Z_l}=\{c(N)\cap s(Z_l,V)\}_{k-d} + \Bbb A_{Z_l},
\tag 3.8
$$
where $\Bbb A_{Z_l}$ are the adjunct terms given by
$$
\Bbb A_{Z_l} =
\sum_{i=0}^{d-2}\sum_{j=1}^{d-1-i}{d-1-i \choose j}
c_i(N)\cap \pi_{l*}(s_{d-i-j}(E_l)\cap s_{k-j}(R_l,V_l)),
\tag 3.9
$$
for $l$ equal to $1$ and $ 2$,
respectively.
\endproclaim{}

\example{Remarks}
(1) It is not required that $Z_1$ and $Z_2$ are residual schemes to
each other in $W$.
However,
it is easy to see that
the hypothesis of the theorem above does imply that,
for $l$ equal to $1$ and $ 2$,
the residual schemes $R_l$ to $E_l$ in $\pi_l^{-1}(W)$
are equal to the inverse image schemes $\pi_l^{-1}(Z_{\hat l})$,
respectively,
where
$$
\hat 1= 2 \quad \text {and} \quad \hat 2 = 1.
\tag 3.10
$$
The above notation $\hat l$ will be use from time to time in this paper.

(2) Theorem 3.6 allows us to write
down many different formulas for our residual decomposition.
The form given in the statement above
is a mixture of the formulas obtained from
two of the three procedures.
A more natural way to express
the decomposition symmetrically is to use the following
formulas obtained directly using the blow-up $\pi$.
For $l$ equal to $1$ and $ 2$,
$$
\aligned
\Bbb R_{Z_l}&=b_*(c_{d-1}(\pi^*N-\widetilde E_1\otimes \widetilde E_2)
\cap [\widetilde E_l])\\
&=\sum_{i=0}^{d-1}\sum_{j=0}^{d-1-i}{d-1-i \choose j}
c_i(N)\cap \pi_*((-\widetilde E_1)^j(-\widetilde E_2)^{d-1-i-j}
\cap [\widetilde E_l])\\
&=\{c(N)\cap s(Z_l,V)\}_{k-d} + \Bbb A_{Z_l},
\endaligned{}
\tag 3.11
$$
where $\Bbb A_{Z_l}$ are adjunct terms defined by
$$
\aligned
\Bbb A_{Z_l}=\sum_{i=0}^{d-2}\sum_{j=1}^{d-1-i}{d-1-i \choose j}
c_i(N)\cap \pi_{*}((-\widetilde E_{\hat l})^j(-\widetilde E_l)^{d-1-i-j}
\cap [\widetilde E_l]).
\endaligned{}
\tag 3.12
$$
Since they are somehow derived directly from the definition
of $X \cdot V$,
they might be more useful in some applications.
\endexample{}

\demo{Proof of Theorem 3.6}
As defined in Theorem 2.4,
let
$$
X \cdot V_1' = \Bbb R_1 + \Bbb R_1'
$$
be the natural decomposition of
$X \cdot V_1'$
with respect to
$E_1'$ and $R_1'$,
$$
X \cdot \widetilde V = \widetilde{\Bbb R}_1 + \widetilde{\Bbb R}_2
$$
be the natural decomposition of
$X \cdot \widetilde V$
with respect to
$\widetilde E_1$ and $\widetilde E_2$,
and
$$
X \cdot V_2' = \Bbb R_2 + \Bbb R_2'
$$
be the natural decomposition of
$X \cdot V_2'$
with respect to
$E_1'$ and $R_2'$.
Since
$$
X \cdot V=
(\pi_1b_1)_*(X \cdot V_1')
=\pi_*(X\cdot \widetilde V)
=(\pi_2b_2)_*(X \cdot V_2')
$$
and
$$
p_*(X \cdot V_1')
=X\cdot \widetilde V
=q_*(X \cdot V_2'),
$$
it is enough to show that
$$
p_*(\Bbb R_1)=\widetilde{\Bbb R}_1
\quad \text {and} \quad
q_*(\Bbb R_2)=\widetilde{\Bbb R}_2.
$$
By the symmetric,
it is enough to show one of them as follows.
$$
\aligned
p_*(\Bbb R_1)
&=p_*(\{c((\pi_1b_1)^*N - E_1'\otimes R_1')
\cap [E_1']\}_{k-d})\\
&=p_*(\{c((\pi p)^*N-p^*(\widetilde E_1\otimes \widetilde E_2))
\cap [E_1']\}_{k-d})\\
&=p_*(\{c(p^*(\pi^*N-\widetilde E_1\otimes \widetilde E_2))
\cap [E_1']\}_{k-d})\\
&=\{c(\pi^*N-\widetilde E_1\otimes \widetilde E_2)
\cap p_*([E_1'])\}_{k-d}\\
&=\{c(\pi^*N-\widetilde E_1\otimes \widetilde E_2)
\cap [\widetilde E_1]\}_{k-d}\\
&=\widetilde{\Bbb R}_1.
\endaligned{}
$$
\qed
\enddemo{}

The hypothesis of Theorem 3.4 (or Theorem 3.6)
is not unreasonable and is often satisfied
in many applications.
Given a subscheme $Z_1$ of $W$,
we are free to pick up a $Z_2$ such that the hypothesis are
satisfied,
since they do not have to be the residual schemes
to each other.
For instance,
for $Z_1$ equal to $p$ or $R(p)$ in Example 3.1,
we can take $Z_2$ to be $p$ or $\emptyset$,
respectively.
Theorem 3.6 can then be applied.
Another useful case is when
$Z_1$ and $Z_2$ are complementary components in $W$.
Roughly speaking,
if a component $Z_1$ of $W$ is given,
we then take $Z_2$ to be the union of all irreducible components of $W$
not appeared in $Z_1$.
It is easy to see that
$Z_1$ and $Z_2$ do form the
residual subschemes with respect to each other in $W$
in this case.

As an immediate consequence of the theorem,
we see that the new decomposition
behaves well with respect to the canonical decomposition of $X\cdot V$.
For example,
the following corollary will be useful for our applications
in Section 4.

\proclaim{Corollary 3.13}
If $Z_1$ and $Z_2$ are irreducible components of $W$ and
there is no contribution to $X\cdot V$ from $Z_1\cap Z_2$,
then the refined residual intersection decomposition coincides
with the canonical decomposition.
\endproclaim{}

Another interesting case is when dimensions
of intersections of distinguished varieties are small.
Since the adjunct term in
each $\Bbb R_{Z_l}$ is a class supported on
the intersection of $Z_1$ and $Z_2$,
that complicated part
will be zero if the dimension of the intersection is less than $k-d$.
In this case,
only $\Bbb M_{Z_l}$ are left and
our decomposition becomes very simple.
In other words,
when the dimension of the intersection of two
components of a intersection product
is less than the dimension of the product itself,
we may treat them as they were connected components disjoint to each other.
In general,
the following corollary follows easily by the induction.

\proclaim{Corollary 3.14}
Let $X\cdot V$ be the refined intersection product defined
from the fiber square (2.1)
$$
\CD
W @>>> V\\
@VfVV @VVfV\\
X @>>i> Y
\endCD{}
$$
and $W_i$'s be irreducible components of $W$.
If for any pair of components $W_i$ and $W_j$,
$i\not = j$,
$$
\codim(W_i\cap W_j,V) > \codim(X,Y)=d,
$$
then the canonical decomposition of $X\cdot V$ is given by
$$X\cdot V=
\sum_i\{c(N)\cap s(W_i,V)\}_{k-d}.
\tag 3.15
$$
\endproclaim{}

\example{Remarks}
(1) In general,
$X\cdot V$ may have many additional distinguished varieties
other than $W_i$'s.
Just assuming all contributions not from $W_i$'s to be zero
is not sufficient for formula (3.15) to hold
since that dose not imply the vanishing of
the adjunct terms as we will see
from examples in Section 4.

(2) The formula (3.15) also follows easily from the fiber diagram
below:
$$
\CD
W' @>>> V'\\
@VpVV @VVpV\\
W @>>> V\\
@VfVV @VVfV\\
X @>>i> Y,
\endCD{}
$$
where
$$
V' = V - \cup_{i \not = j} (W_i\cap W_j),
\quad W' = W - \cup_{i \not = j} (W_i\cap W_j),
$$
and $p$ is the open embedding.
\endexample{}

If both $Z_1$ and $Z_2$ are regularly embedded in $V$,
then we can further reduce our formulas
into a more explicit form that uses no push-forwards.
Those formulas are useful in actual computations,
since they involve mostly the characteristic classes of vector bundles.

\proclaim{Corollary 3.16}
In the set-up of Theorem 3.6,
assume further that $Z_1$ and $Z_2$ intersect properly (transversely).
If both $Z_1$ and $Z_2$ are regularly embedded in $V$
of codimensions $r_1$ and $r_2$
with the normal bundles $N_1$ and $N_2$,
respectively,
then we have the following new formulas for
$\Bbb M_{Z_l}$ and $\Bbb A_{Z_l}$.
For $l$ equal to $1$ and $2$,
$$
\Bbb M_{Z_l} = \{c(N)\cap s(Z_l,V)\}_{k-d}
=c_{d-r_l}(N-N_l)\cap [Z_l],
\tag 3.17
$$
$$
\Bbb A_{Z_l} =-\sum_{i=0}^{d-r_1-r_2}\sum_{j=r_{\hat l}}^{d-r_l-i}
{d-1-i \choose j}
c_i(N)s_{j-r_{\hat l}}(N_{\hat l})s_{d-r_l-i-j}(N_l)\cap [Z_{int}],
\tag 3.18
$$
where $Z_{int}$ is the intersection of $Z_1$ and $Z_2$.
\endproclaim{}

\demo{Proof}
By the symmetry,
it is enough to verify the new formulas for $\Bbb M_{Z_1}$
and $\Bbb A_{Z_1}$.
The formula (3.17) for $\Bbb M_{Z_1}$ is easy to see,
since
$$
s(Z_1,V) = s(N_1)\cap [Z_1] = c^{-1}(N_1)\cap [Z_1].
$$
To see the new formula for the adjunct term $\Bbb A_{Z_1}$,
we will start from its formula given in Theorem 3.6.
Notice from (3.9) that $\pi_{1*}$ pushes a class supported on $E_1$
into a class supported on $Z_{int}$.
We hence may consider $\pi_1$ as the projection from
the projective bundle $E_1=\Bbb P(N_1)$
to $Z_1$.
We claim that
$$
E_1\cap s(R_1, V_1)=\pi_1^*(s(N_2)\cap[Z_{int}]).
\tag 3.19
$$
To see how to get the formula from the claim above,
notice that $Z_{int}$ has codimension $r_1+r_2$ in $V$.
We will assume that it is less the $d$,
since $\Bbb A_{Z_1}$ is equal to zero otherwise.
We thus have from (3.9) that
$$
\aligned
\Bbb A_{Z_1}&=\sum_{i=0}^{d-2}\sum_{j=1}^{d-1-i}{d-1-i \choose j}
c_i(N)\cap \pi_{1*}(s_{d-i-j}(E_1)\cap s_{k-j}(R_1,V_1))\\
&=-\sum_{i=0}^{d-2}\sum_{j=1}^{d-1-i}{d-1-i \choose j}
c_i(N)\cap \pi_{1*}((-E_1)^{d-i-j-1}\cap (E_1\cap s(R_1,V_1))_{k-j-1})\\
&=-\sum_{i=0}^{d-2}\sum_{j=1}^{d-1-i}{d-1-i \choose j}
c_i(N)\cap \pi_{1*}((-E_1)^{d-i-j-1}\cap
\pi_{1}^*(s_{j-r_2}(N_2)\cap[Z_{int}]))\\
&=-\sum_{i=0}^{d-2}\sum_{j=1}^{d-1-i}{d-1-i \choose j}
c_i(N)s_{d-r_1-i-j}(N_1)s_{j-r_2}(N_2)\cap [Z_{int}].
\endaligned{}
$$
Dropping the terms which contain any factor with a negative index,
we then obtain the formula (3.18).
We hence need only to verify the claim (3.19).
To see this,
notice from the hypotheses that
$$
\pi_1^*(s(N_2)\cap[Z_{int}])
=\pi_1^*(s(Z_{int},Z_1))
=s(R_1\cap E_1, E_1).
\tag 3.20
$$
For any pair $X$ in $Y$,
let $N(X,Y)$ be
the normal cone to
$X$ in $Y$.
Restricting everything to $R_1\cap E_1$
and notice that both $N(E_1, V_1)$ and $N(R_1\cap E_1, R_1)$
are isomorphic to $E_1$,
we have the following exact sequences
$$
0 \to N(R_1\cap E_1, E_1)
\to N(R_1\cap E_1, V_1)
\to E_1 \to 0
$$
and
$$
0 \to E_1 \to N(R_1\cap E_1, V_1)
\to N(R_1, V_1)
\to 0.
$$
Therefore,
continuing from (3.20),
we get
$$
\aligned
s(R_1\cap E_1, E_1)
&=s(N(R_1\cap E_1, E_1))\\
&=c(E_1)\cap s(N(R_1\cap E_1, V_1))\\
&=c(E_1)\cap(s(E_1)\cap s(N(R_1, V_1)))\\
&=E_1\cap s(R_1, V_1).
\endaligned{}
$$
This completes the proof.
\qed
\enddemo{}

\smallpagebreak
\heading 4. Geometric meaning and applications ---
degenerations of hypersurfaces and their limiting $\Bbb P^r$
\endheading{}
\smallpagebreak

An important reason that inspired our study of the new
residual intersection is its application to some geometric problems.
As we have mentioned in the introduction,
using dynamic intersections,
we can interpret the refined intersection product $X\cdot V$
as a limit derived from deformation theory.
In this point of view,
our decomposition has nice geometric meaning
and tells us how the limit
is distributed in different components.
Instead of doing it in more general terms,
we will study the case of the limiting $\Bbb P^r$'s
in hypersurfaces of $\Bbb P^n$.
The same method can be applied to many other settings.

We will first recall some definitions and facts from [W2].
Let $X$ be a hypersurface of degree $d$ in $\Bbb P^n$
and $F_X$ be the scheme of $\Bbb P^r$'s contained in $X$.
If $X$ is generic,
then $F_X$ will have the expected dimension (or empty)
and its class in the Chow group of $G$
is given by the top Chern class of the vector bundle ${\Sym}^dU^*$,
where $G$ is the Grassmannian of $\Bbb P^r$'s in $\Bbb P^n$
and $U$ is the universal subbundle on $G$.
When we deform hypersurfaces into a degenerate $X$,
the dimension of $F_X$ can jump.
In this case,
there is a subscheme $F_{lim}$ of $F_{X}$
with the expected dimension
which consists of the limiting
$\Bbb P^r$'s in $X$ with respect to a general deformation.
In general,
the structure of $F_{lim}$ is much more complicated
than the structure of $F_X$.
However,
to find the class of $F_{lim}$ and its distribution
in different components of $F_{X}$,
one needs only (at least in theory)
to understand the structure of $F_X$.
To see this,
we consider the following fiber square:
$$
\CD
F_X @>>> G\\
@VVV @VVs_XV\\
G @>>i> {\Sym}^dU^*
\endCD{}
$$
where
$i$ is the zero-section embedding and
$s_X$ is the section of ${\Sym}^dU^*$ induced by $X$.
As it turns out,
the class of $F_{lim}$ is equal to the refined intersection
product $G\cdot G$ defined from the fiber square above
and the distribution of the limiting $\Bbb P^r$'s
in different
components of $F_{X}$ is given by the canonical decomposition of
$G\cdot G$.
For details,
please refer to [W2]
($G\cdot G$ is called $\Cal R_X$ there.)
Therefore,
the problem is reduced into finding such a decomposition.
For the obvious reason,
we did not carry out the program along this line in [W2].
Instead,
we found the decomposition by conducting
a direct study of $F_{lim}$ in the given cases.
Unfortunately,
the methods used there are not easy to be generalized.
However,
it did give some hints about the possibility of a general formula
and this paper is a result of the search that follows.
With the results in the previous sections,
we are ready to apply the
refined residual intersection formulas
to compute the distribution of limiting $\Bbb P^r$'s in $F_X$.
The point is that,
since our formulas use only information
on $F_X$,
we now have a general way to find the distribution
of the limits without doing any study about $F_{lim}$ itself.
To simplify the notation,
we will omit $\cap [G]$ in our formulas to
identify an element in $A^*(G)$ with its dual in
$A_*(G)$.
Furthermore,
for any vector bundle $E$ and any integer $m$,
we will use $E(m)$ to denote a vector bundle of the same rank
such that its Chern class is given by
$$
c(E(m))=\Adams(m,c(E))=\sum_i m^ic_i(E).
$$
where $\Adams(*,*)$ is the Adams operator of $K$-theorem.
Notice that the Segre class of $E(m)$ is thus determined by
the Segre class of $E$ via
$$
s(E(m))=\Adams(m,s(E))=\sum_i m^is_i(E).
$$
For any positive integer $m$,
we will also use the following notation
$$
r_m=\rank({\Sym}^mU^*)={m+r \choose r}.
$$

\proclaim{Proposition 4.1}
Let $X_k^e$ be the $e$-fold of a generic hypersurface of degree $k$
in $\Bbb P^n$
and $X_l^f$ be the $f$-fold of a generic hypersurface of degree $l$
in $\Bbb P^n$.
If we deform a generic hypersurface of degree $d$
in $\Bbb P^n$ into the union of
$X_k^e$ and $X_l^f$ (with $ke+lf=d$),
then the class $[F_{lim}(X_k^e)]$
of the limiting $\Bbb P^r$'s in $X_k^e$
can be computed by the following formulas
$$
[F_{lim}(X_k^e)]= \Bbb R_{F_{X_k^e}}
= \Bbb M_{F_{X_k^e}} + \Bbb A_{F_{X_k^e}},
\tag 4.2
$$
where
$$
\aligned
\Bbb M_{F_{X_k^e}} = & c_{r_k}({\Sym}^kU^*(e))
c_{r_d-r_k}({\Sym}^dU^*-{\Sym}^kU^*(e))\\
= & c_{r_k}({\Sym}^kU^*(e))
\sum_{i=0}^{r_d-r_k}c_i({\Sym}^dU^*)c_{r_d-r_k-i}({\Sym}^kU^*(e)),
\endaligned{}
\tag 4.3
$$
and
$$
\aligned
\Bbb A_{F_{X_k^e}}= & -c_{r_k}({\Sym}^kU^*(e))c_{r_l}({\Sym}^lU^*(f))
\sum_{i=0}^{r_d-r_k-r_l}\sum_{j=r_l}^{r_d-r_k-i}{r_d-1-i \choose j}\\
& \times c_{i}({\Sym}^dU^*)
s_{j-r_l}({\Sym}^lU^*(f))s_{r_d-r_k-i-j}({\Sym}^kU^*(e)).
\endaligned{}
\tag 4.4
$$
Similarly,
the class $[F_{lim}(X_l^f)]$
of the limiting $\Bbb P^r$'s in $X_l^f$
is given by
$$
[F_{lim}(X_l^f)]= \Bbb R_{F_{X_l^f}}
= \Bbb M_{F_{X_l^f}} + \Bbb A_{F_{X_l^f}},
\tag 4.5
$$
where
$$
\aligned
\Bbb M_{F_{X_l^f}}
= & c_{r_l}({\Sym}^lU^*(f))c_{r_d-r_l}({\Sym}^dU^*-{\Sym}^lU^*(f))\\
= & c_{r_l}({\Sym}^lU^*(f))
\sum_{i=0}^{r_d-r_l}c_{i}({\Sym}^dU^*)c_{r_d-r_l-i}({\Sym}^lU^*(f))
\endaligned{}
\tag 4.6
$$
and
$$
\aligned
\Bbb A_{F_{X_l^f}} = & - c_{r_l}({\Sym}^lU^*(f))c_{r_k}({\Sym}^kU^*(e))
\sum_{i=0}^{r_d-r_k-r_l}\sum_{j=r_k}^{r_d-r_l-i}{r_d-1-i \choose j}\\
& \times c_{i}({\Sym}^dU^*)
s_{j-r_k}({\Sym}^kU^*(e))s_{r_d-r_l-i-j}({\Sym}^lU^*(f)).
\endaligned{}
\tag 4.7
$$
\endproclaim{}

\demo{Proof}
Since $X$ is the union of $X_k^e$ and $X_l^f$,
$F_X$ consists of two components $F_{X_k^e}$ and $F_{X_l^f}$.
We can therefore apply
the refined residual intersection decomposition to
$$
G\cdot G = [F_{lim}]
$$
with respect to the decomposition
$$
F_X= F_{X_k^e} \cup F_{X_l^f}
$$
of $F_X$.
As sets,
$F_{X_k^e}$ and $F_{X_l^f}$ are equal to
$F_{X_k}$ and $F_{X_l}$,
respectively.
Furthermore,
the Segre classes of $F_{X_k^e}$ and $F_{X_l^f}$
in $G$ is closely related to the normal
bundles of $F_{X_k}$ and $F_{X_l}$ in G,
respectively.
In fact,
by the same method used in the proof of Theorem 1.4 of [W3],
we can easily get the following:
$$
s(F_{X_k^e}, G)=\Adams(e,s(F_{X_k},G))=
s({\Sym}^kU^*(e))\cap [F_{X_k^e}]
$$
and
$$
s(F_{X_l^f}, G)=\Adams(f,s(F_{X_l},G))=
s({\Sym}^lU^*(f))\cap [F_{X_l^f}].
$$
Therefore,
$F_{X_k^e}$ and $F_{X_l^f}$ can be considered as
being regularly embedded in $G$ with normal
bundles
${\Sym}^kU^*(e)$ and ${\Sym}^lU^*(f)$,
respectively.
Since $X_k^e$ and $X_l^f$ are generic,
$F_{X_k^e}$ and $F_{X_l^f}$ do intersect properly
and the class of their intersection is given by
$$
\aligned
[F_{X_k^e} \cap F_{X_l^f}]&=
c_{top}({\Sym}^kU^*(e)\oplus {\Sym}^lU^*(f))\cap [G]\\
&=c_{top}({\Sym}^kU^*(e))c_{top}({\Sym}^lU^*(f)) \cap [G].
\endaligned{}
$$
Now,
apply Corollary 3.16.
\qed
\enddemo{}

Proposition 4.1 generalizes
some earlier results on the
limiting $\Bbb P^r$'s from [W1], [W2], and [W3].
Notice that the corresponding formulas in [W1], [W2], and [W3]
are obtained from very different methods.

(a) If $r=1$ and $e=f=1$,
then the adjunct parts are equal to zero by the dimension count.
Our formulas hence are reduced to and coincide with the ones given in [W1]
in the case of the limiting lines.

(b) More general,
if just set $e=f=1$,
we then have the case studied in [W2].
Notice that our formulas in this case are different from (and simpler than)
the corresponding formulas in [W2].
Of course,
they must be equal to each other but that
does not seem to be obvious if one just looks at the formulas
themselves.

(c) If $f=0$,
we then have the case studied in [W3] and
our formulas is thus a generalization of the formulas obtained there.

(d) The set-ups and the formulas are new if $ef > 1$.
Notice that geometric and infinitesimal methods used in [W1] and [W2]
will be difficult to apply in such cases,
since we are dealing with non-reduced scheme structures
and $F_{lim}$
will be very hard to understand.

The formulas in Proposition 4.1 are not as complicated
as their look and can be calculated with routine procedures
in intersection theory.
Our examples below
will be mainly for the new case (d) and the case (b)
using the new formulas.
We will write $X_k$ for $X_k^1$.
Most of the calculations are done on computers
using the Maple package
schubert [KS] written by Katz and Str{\o}mme.
It is quite easy to write schubert code for our formulas
and we will be happy to send the copies upon request.
But first,
we will check a simple example by hand.

\proclaim {Example 4.8} Degenerations of
cubic hypersurfaces in $\Bbb P^n$, $n\geq 3$, and their lines
\endproclaim{}
The Chow ring of $G$ of lines in $\Bbb P^n$ is
generated by $c_1(U^*)$ and $c_2(U^*)$.
Let
$$
c_1(U^*)=x \quad \text {and} \quad
c_2(U^*)=y.
$$
It is an easy computation by the splitting principle that
$$
\left\{
\aligned
c_1({\Sym}^3U^*)&=6x,\\
c_2({\Sym}^3U^*)&=11x^2+10y,\\
c_3({\Sym}^3U^*)&=6x^3+30xy,\\
c_4({\Sym}^3U^*)&=18x^2y+9y^2.
\endaligned{}
\right .
\tag 4.9
$$
and
$$
s_1(U^*)=-x,\quad
s_2(U^*)=x^2-y,\quad
\cdots.
\tag 4.10
$$
Therefore,
the class of lines on a generic
cubic $X_3$ is equal to
$$
[F_{X_3}]=c_4({\Sym}^3U^*)=18x^2y+9y^2.
$$
To compute the distribution of the limiting lines
when we deform $X_3$ into the union
of a $(n-1)$-plane $X_1$ and a double-$(n-1)$-plane $X_1^2$,
we set
$$
r=1,\quad d=3, \quad l=k=e=1, \quad \text {and} \quad f=2
$$
in the formulas in Proposition 4.1.
This together with (4.9) and (4.10) give us
$$
\aligned
[F_{lim}(X_1)]=& \Bbb M_{F_{X_1}} + \Bbb A_{F_{X_1}}\\
=&c_{2}(U^*)c_{4-2}({\Sym}^3U^*-U^*) -
c_{2}(U^*)c_{2}(U^*(2)){4-1 \choose 2}\\
=&c_2(U^*)\sum_{i=0}^2c_{i}({\Sym}^3U^*)s_{2-i}(U^*)
-12c_{2}(U^*)c_{2}(U^*)\\
=&y(x^2-y+6x(-x)+11x^2+10y)-12y^2\\
=&6x^2y-3y^2
\endaligned{}
\tag 4.11
$$
and
$$
\aligned
[F_{lim}(X_1^2)]=& \Bbb M_{F_{X_1^2}} + \Bbb A_{F_{X_1^2}}\\
=&c_{2}(U^*(2))c_{4-2}({\Sym}^3U^*-U^*(2)) -
c_{2}(U^*)c_{2}(U^*(2)){4-1 \choose 2}\\
=&c_2(U^*(2))\sum_{i=0}^2c_{i}({\Sym}^3U^*)s_{2-i}(U^*(2))
-12c_{2}(U^*)c_{2}(U^*)\\
=&4y(4(x^2-y)+6x(-2x)+11x^2+10y)-12y^2\\
=&12x^2y+12y^2.
\endaligned{}
\tag 4.12
$$
Notice that we do have
$$
[F_{lim}(X_1)]+[F_{lim}(X_1^2)]=18x^2y+9y^2=[F_{X_3}]
$$
as expected.

In particular,
if $n=3$,
then
$$
x^2y=y^2=[point].
$$
Therefore,
we see that $24$ of the $27$ lines in a generic cubic surface
go to the double-plane $X_1^2$ and other $3$ lines go to the plane
$X_1$ as one takes the limit.
In fact,
it is easy to identify those lines geometrically.
Let
$$
C=X_1\cap D \subset X_1
$$
be the elliptic curve in the plane $X_1$,
where $D$ is a cubic surface determined
infinitesimally by the degeneration.
It intersects $X_1^2$ at three points.
Three limiting lines in $X_1$ are exactly the tangent lines
in $X_1$ to $C$ at those three points.
Similarly,
we have $3$ tangent lines in $X_1^2$.
That gives $24$ limiting lines in $X_1^2$
since now each lines has the multiplicity of
$8$ on account of the non-reduced
scheme structure of $X_1^2$.

Another simple example that can be easily checked by hand is
the case of the limiting $\Bbb P^2$'s in a degeneration of quadrics.
Interested readers may try this example as a comparison
to the computations made in [W2] [Example 3, Section 5]
in which more complicated formulas are used.

\proclaim {Example 4.13} Degenerations of
quintic threefolds and their lines.
\endproclaim{}

It is well-known that there are $2875$
lines on a generic
generic quintic threefold.
There are seven different
degenerations of a generic $X_5$ into $X_k^e \cup X_l^f$.
Two of such cases are studied
in [K] and [W1] so we will only list below the results
for five new cases with
at least one non-reduced component.
It will be interesting to have those numbers checked
geometrically.
To compare with
the standard residual intersection formula,
we will also list values for both
the main terms and the adjunct terms.

\roster
\itemitem{} Case 1, $X_5 \to X_1^4 \cup X_1$.
$$
[F_{lim}(X_1^4)] =
\Bbb R_{F_{X_1^4}}=\Bbb M_{F_{X_1^4}} + \Bbb A_{F_{X_1^4}}
= 2400 + 320 = 2720;
$$
$$
[F_{lim}(X_1)] =
\Bbb R_{F_{X_1}}=\Bbb M_{F_{X_1}} + \Bbb A_{F_{X_1}}
= 1275 - 1120 = 155.
$$

\itemitem{} Case 2, $X_5 \to X_1^3 \cup X_2$.
$$
[F_{lim}(X_1^3)] =
\Bbb R_{F_{X_1^3}}=\Bbb M_{F_{X_1^3}} + \Bbb A_{F_{X_1^3}}
= 3195 - 540 = 2655;
$$
$$
[F_{lim}(X_2)] =
\Bbb R_{F_{X_2}}=\Bbb M_{F_{X_2}} + \Bbb A_{F_{X_2}}
= 1300 - 1080 = 220.
$$

\itemitem{} Case 3, $X_5 \to X_1^3 \cup X_1^2$.
$$
[F_{lim}(X_1^3)] =
\Bbb R_{F_{X_1^3}}=\Bbb M_{F_{X_1^3}} + \Bbb A_{F_{X_1^3}}
= 3195 - 1080 = 2115;
$$
$$
[F_{lim}(X_1^2)] =
\Bbb R_{F_{X_1^2}}=\Bbb M_{F_{X_1^2}} + \Bbb A_{F_{X_1^2}}
= 2920 - 2160 = 760.
$$

\itemitem{} Case 4, $X_5 \to X_1^2 \cup X_3$.
$$
[F_{lim}(X_1^2)] =
\Bbb R_{F_{X_1^2}}=\Bbb M_{F_{X_1^2}} + \Bbb A_{F_{X_1^2}}
= 2920 - 540 = 2380;
$$
$$
[F_{lim}(X_3)] =
\Bbb R_{F_{X_3}}=\Bbb M_{F_{X_3}} + \Bbb A_{F_{X_3}}
= 1575 - 1080 = 495.
$$

\itemitem{} Case 5, $X_5 \to X_2^2 \cup X_1$.
$$
[F_{lim}(X_2^2)] =
\Bbb R_{F_{X_2^2}}=\Bbb M_{F_{X_2^2}} + \Bbb A_{F_{X_2^2}}
= 2880 - 640 = 2240;
$$
$$
[F_{lim}(X_1)] =
\Bbb R_{F_{X_1}}=\Bbb M_{F_{X_1}} + \Bbb A_{F_{X_1}}
= 1275 - 640 = 635.
$$
\endroster{}
Notice that in each of the cases above we do have that
$$
[F_{lim}(X_k^e)]+[F_{lim}(X_l^f)]=
\Bbb R_{F_{X_k^e}} + \Bbb R_{F_{X_l^f}} = 2875
$$
as expected.

\proclaim {Example 4.14} Degenerations of
quartic hypersurfaces in $\Bbb P^7$ and their $\Bbb P^2$'s.
\endproclaim{}

As we have computed in [W2],
there are $3,297,280$
$\Bbb P^2$'s on a generic
quartic hypersurface in $\Bbb P^7$.
There are five different degenerations for which Proposition 4.1
can be applied.
We will list results of our computation for all five cases.
Notice that the values of the main terms can range
from negative ones to those that are bigger than the whole class.
In particular,
we see that the adjunct terms are not equal to zero
even though there is no
contribution from the variety on which they are supported (see [W2]).

\roster
\itemitem{} Case 1, $X_4 \to X_3 \cup X_1$.
$$
[F_{lim}(X_3)] =
\Bbb R_{F_{X_3}} = \Bbb M_{F_{X_3}} + \Bbb A_{F_{X_3}}
=3,304,098 - 2,820,258 = 483,840;
$$
$$
[F_{lim}(X_1)] =
\Bbb R_{F_{X_1}}=\Bbb M_{F_{X_1}} + \Bbb A_{F_{X_1}}
=3,656,569 - 843,129 = 2,813,440.
$$

\itemitem{} Case 2, $X_4 \to X_2 \cup X_2$.
$$
[F_{lim}(X_2)] =
\Bbb R_{F_{X_2}}=\Bbb M_{F_{X_2}} + \Bbb A_{F_{X_2}}
=3,087,616 - 1,438,976 = 1,648,640.
$$

\itemitem{} Case 3, $X_4 \to X_1^3 \cup X_1$.
$$
[F_{lim}(X_1^3)] =
\Bbb R_{F_{X_1^3}}=\Bbb M_{F_{X_1^3}} + \Bbb A_{F_{X_1^3}}
=-20,855,205 + 24,000,165 = 3,144,960;
$$
$$
[F_{lim}(X_1)] =
\Bbb R_{F_{X_1}}=\Bbb M_{F_{X_1}} + \Bbb A_{F_{X_1}}
=3,656,569 - 3,504,249 = 152,320.
$$

\itemitem{} Case 4, $X_4 \to X_1^2 \cup X_1^2$.
$$
[F_{lim}(X_1^2)] =
\Bbb R_{F_{X_1^2}}=\Bbb M_{F_{X_1^2}} + \Bbb A_{F_{X_1^2}}
=2,645,888 - 997,248 =1,648,640.
$$

\itemitem{} Case 5, $X_4 \to X_1^2 \cup X_2$.
$$
[F_{lim}(X_1^2)] =
\Bbb R_{F_{X_1^2}}=\Bbb M_{F_{X_1^2}} + \Bbb A_{F_{X_1^2}}
=2,645,888 + 561,792 =3,207,680;
$$
$$
[F_{lim}(X_2)] =
\Bbb R_{F_{X_2}}=\Bbb M_{F_{X_2}} + \Bbb A_{F_{X_2}}
=3,087,616 - 2,998,016 = 89,600.
$$
\endroster{}
Notice that
the results in case (1) and case (2) do coincide with the ones given
in Example 2 of [W2]
obtained from different methods and formulas.

We have calculated many other examples.
In the cases that other methods and formulas are available,
such as the case of Example 4 of [W2] as interested readers may check,
the results do coincide as expected.

Our last application is about a family of identities for
the characteristic classes of vector bundles.
For any vector bundle $E$,
we will formally extend the definition of
$s_i(E)$ to negative $i$ by
$$
s_i(E)=0, \quad -\rank (E)< i < 0, \quad \quad
s_i(E)=- 1/c_{top}(E), \quad i=-\rank (E).
$$

\proclaim{Corollary 4.15}
Let $E$ be a vector bundle of rank $r+1$ on a purely dimensional
scheme $X$.
If $E$ is generated by its sections,
then the following family of identities holds for any set of
non-negative integers $k$, $l$, and $d$
with $d=k+l$.
$$
\aligned
c_{r_d}({\Sym}^dE) = &
-c_{r_k}({\Sym}^kE)c_{r_l}({\Sym}^lE)\times\\
(4.16)\quad \quad
&\left\{ \sum_{i=0}^{r_d-r_k}\sum_{j=0}^{r_d-r_k-i}
{r_d-1-i \choose j}c_{i}({\Sym}^dE)
s_{j-r_l}({\Sym}^lE)s_{r_d-r_k-i-j}({\Sym}^kE)\right .\\
&+ \left .\sum_{i=0}^{r_d-r_l}\sum_{j=0}^{r_d-r_l-i}
{r_d-1-i \choose j}c_{i}({\Sym}^dE)
s_{j-r_k}({\Sym}^kE)s_{r_d-r_l-i-j}({\Sym}^lE)\right \}.
\endaligned{}
$$
\endproclaim{}

\example{Remark}
It is pointed out to the author by the referee that
Corollary 4.15 actually holds for any vector bundle $E$ of rank $r+1$:
Taking an ample line bundle $H$ and replacing
$E$ by $E\otimes H^t$,
we can then consider both sides of (4.16) as a polynomial
in $t$.
Since they agree for $t$ sufficiently large,
they must agree for all $t$.
In particular, they agree for $t=0$.
\endexample{}

\demo{Proof of Corollary 4.15}
Consider the following fiber square
$$
\CD
W @>>> X\\
@VVV @VVfV\\
X @>>i> {\Sym}^dE
\endCD{}
$$
where
$i$ is the zero-section embedding and
$f$ is any section of ${\Sym}^dE$.
Therefore,
$W$ is the zero-scheme of $f$.
Since $d=k+l$,
a section $s_k$ of ${\Sym}^kE$ and
a section $s_l$ of ${\Sym}^lE$ induce a section of ${\Sym}^dE$.
Let $f$ be such a section and we hence have
$$
W=Z_1 \cup Z_2,
$$
where $Z_1$ and $Z_2$ are zero-schemes of $s_k$ and $s_l$,
respectively.
Since $E$ is generated by its sections,
if we take $s_k$ and $s_l$ to be generic then
$Z_1$ and $Z_2$ are regularly embedded in $X$
with normal bundle ${\Sym}^kE$ and ${\Sym}^lE$,
respectively.
Notice the normal bundle of $X$ in ${\Sym}^dE$
is ${\Sym}^dE$ itself and it is well-known that
the image of $X\cdot X$ in the Chow group of $X$ is equal to
$c_{top}({\Sym}^dE)\cap [X]$.
{}From this,
the family of identities follows
directly from our residual intersection formulas.
\qed
\enddemo{}

{}

\enddocument{}

\bye